\documentclass[12pt, a4paper]{amsart}

\usepackage{amssymb, amsthm, amsmath, amsfonts}

\usepackage[top=2.6cm, bottom=2.5cm, left=2.6cm, right=2.6cm]{geometry}

\theoremstyle{plain}
\newtheorem{thm}{Theorem}
\theoremstyle{remark}
\newtheorem{remark}{Remark}

\begin{document}

\title{Non-null framed bordant simple Lie groups}

\begin{abstract}
Let $G$ be a compact simple Lie group equipped with the left invariant framing 
$L$. It is known that there are several groups $G$ such that $(G, L)$ is non-null 
framed bordant. Previously we gave an alternative proof of  these results using  the decomposition formula of its bordism class into a Kronecker product by E. Ossa. In this note we propose a verification formula by reconsidering it, through 
a little more ingenious in the use of this product formula, and try to apply it to 
the non-null bordantness results above. 
\end{abstract}

\author{Haruo Minami}
\address{H. Minami: Professor Emeritus, Nara University of Education}
\email{hminami@camel.plala.or.jp}

\maketitle

\section{Introduction and statement of results}

Let $[G, L]$ be the framed bordism class of a compact simple Lie group $G$ 
equipped with the left invariant framing $L$. There are several well-known 
results on the nonzeroness of $[G, L]$, stated in Theorem 2 below. In this 
note we present a formulation of the method used in ~\cite{M} to verify these  
nonzeroness results, thereby trying again to prove those results. 

Let $S\subset G$ be a circle subgroup. Let $H\subset G$ be a closed 
connected subgroup with $S\subset H$ as a subgroup such that   
$G/H$ is diffeomorphic to the unit sphere $S^m\subset\mathbb{R}^{m+1}$. Identifying $G/H=S^m$ we regard $\pi : G\to G/H$ as the principal 
$H$-bundle over $S^m$ where $\pi$ is the quotient map. Let 
$x=(t_0, \cdots, t_m)\in S^m$ 
where we assume that $\pi(e)=(0, \cdots, 0, 1)$, $e$ denoting the identity 
element of $G$. Let $S^{m-1}$ be the equator of $S^m$ defined by $t_m=0$ 
and let $T : S^{m-1} \to H$ denote the characteristic map of the 
bundle $\pi : G \to S^m$. We assume that $T$ is invariant with 
respect to the involusion 
\[
\lambda :  (t_0, \cdots, t_{m-1}, 0)\to  (-t_0, \cdots, -t_{m-1}, 0).
\]
Suppose given a closed connected subgroup $K\subset H$ such that   
$S\subset K$ is a subgroup and $H/K$ is diffeomorphic to a sphere $S^{r-1}$ 
for some $2\le r \le m$. Similarly to the case above we identify $H/K=S^{r-1}$ 
and regard $p : H\to H/K$ as the principal $K$-bundle over $S^{r-1}$ where $p$ denotes the quotient map. Consider the composite 
$p\circ T : S^{m-1}\to H\to S^{r-1}$. Then its homotopy class can be expressed 
as a multiple of a generator $\alpha\in\pi_{m-1}(S^{r-1})$, i.e. 
$[p\circ T]=d\alpha$ 
for some $0\le d\le \mathrm{ord}(\alpha)-1$.  
\begin{thm}
Suppose $\pi_{m-1}^S=0$. Then if $[G, L]=0$, then we have
\begin{equation*}
 (1+(-1)^{m-1}d(q+1))[H, L]=0\qquad \text{for  some} \ \ 
0\le q\le\left|\pi_{m+1}(S^2)\right|-1. \tag{$*$}
\end{equation*}
\end{thm}

Let $M_n=SO(n), SU(n)$ or $Sp(n)$. In the above, if we take $G=M_n$, 
then $H=M_{n-1}$ and $K=M_{n-2}$ where 
$M_{n-i}=M_{n-i}\times\{I_i\}$, $I_i$ being the identity matrix 
of size $i$. From the matrix form on ~\cite[pp.120, 125, 130]{St} 
we know that the characteristic map of 
the principal $M_{n-1}$-bundle of $M_n$ over $S^{n-1}, 
S^{2n-1}$ or $S^{4n-1}$ satisfies the symmetry property above 
(by following the coordinate rule there). In particular, in the case 
$G=SU(n)$ its characteristic map $T : S^{2n-2}\to SU(n-1)$ is given by 
\[T(x)=T'(x)\begin{pmatrix}
I_{n-2} & 0 \\
0 & -(1+z_{n-1})^2(1-z_{n-1})^{-2} \\
\end{pmatrix}
\]
via the homeomorphism between $U(n)$ and $S^1\times SU(n)$ by modifying 
that of $\pi : U(n)\to S^{2n-1}$, denoted by $T'$. In the case $G=G_2$ there is 
a subgroup $SU(3)$ with $G_2/SU(3)=S^6$. This allows us to take $G=SU(3)$ and $K=SU(2)\subset SU(3)$. Then the characteristic map 
$T : S^5\to H=SU(3)$ can be determined by comparing with that of 
$\pi : SO(7)\to S^6$ through the inclusion $SU(3) \to SO(6)$.
 
Using these facts and Theorem 1, based on the results of calculation of 
$\pi_{n+k}(S^n)$ in ~\cite{T} and ~\cite{MT}, we obtain the following theorem. 
\begin{thm}[\!\!~\cite{A}, ~\cite{BS}, ~\cite{G}, ~\cite{K}, ~\cite{S}, ~\cite{W}]{}
\begin{align*} 
& \ \mathrm{(i)} \ \pi_3^S=\mathbb{Z}_{24}[SU(2), L], \\
& \, \mathrm{(ii)} \ \pi_8^S=\mathbb{Z}_2\oplus\mathbb{Z}_2[SU(3), L], \\ 
&\mathrm{(iii)} \ \pi_{10}^S=\mathbb{Z}_2\oplus\mathbb{Z}_3[Sp(2), L], \\
&\mathrm{(iv)} \ \pi_{14}^S=\mathbb{Z}_2\oplus\mathbb{Z}_2[G_2, L], \\
& \, \mathrm{(v)} \ \pi_{15}^S=\mathbb{Z}_{240}\oplus\mathbb{Z}_2[SU(4), L], \quad
[SU(4), L]=\eta [G_2, L],\qquad\qquad\qquad\qquad\qquad\quad \\
&\mathrm{(vi)} \ \pi_{21}^S=\mathbb{Z}_2\oplus\mathbb{Z}_2[Sp(3), L].\\ 
\end{align*} 
\end{thm}

\begin{remark}
There are no other simple Lie 
groups $G$ such that $[G,  L]\ne 0$ except for the ones mentioned in 
Theorem 2 above  ~\cite{M2} (cf.  ~\cite{M3}).
\end{remark}

\section{Proof of Theorem 1}

Let $E_\pm=\{x\in S^m| \pm t_m\ge 0\}$ be the hemispheres of $S^m$ with 
$\partial E_\pm=S^{m-1}$. Let 
$\tau : \partial E_+\times H\to \partial E_-\times H$ be the bundle isomorphism given by $(x, y)\mapsto (x, T(x)y)$,  $x\in \partial E_+$, $y\in H$. By gluing the products bundles $E_\pm\times H$ under this isomorphism we obtain a bundle decomposition of $\pi : G\to S^m$   
\[
G \cong (E_+\times H)\cup_\tau (E_-\times H).
\]
For the associated bundle $\bar{\pi} : G/S\to S^m$ 
we have a similar decomposition 
\[
G/S \cong (E_+\times H/S)\cup_{\bar{\tau}} (E_-\times H/S)
\]
where $\bar{\tau} : \partial E_+\times H/S\to \partial E_-\times H/S$ 
is the bundle isomorphism given by 
$(x, yS)\mapsto (x, T(x)yS)$, $x\in \partial E_+$, $y\in H$.

Consider an embedding $G/S\hookrightarrow \mathbb{R}^{k+n-1}$ 
for sufficiently large $k>0$ where $n=\dim G$. Since $G/S$ can be 
regarded as a framed manifold equipped with a framing induced by $L$ on 
$G$ ~\cite{LS} (cf. ~\cite{K}, ~\cite{O}), the normal bundle $\nu$ of 
$G/S$ in $\mathbb{R}^{k+n-1}$ becomes isomorphic to the product 
bundle $\mathbb{R}^k\times G/S$ over $G/S$. So via the Pontrjagin-Thom 
construction we obtain a collapse map $f : S^{k+n-1}\to S^k(G/S^+)$ which 
gives $[G/S]\in \pi^S_{n-1}(G/S^+)$, the fundamental bordism class of $G/S$. 

Let $\xi$ denote the complex line bundle $G\times_SV\to G/S$ associated to 
the standard complex representation $V$ of 
$S^1=\{z\in\mathbb{C}| |z|=1\}$ via an isomorphism $S\cong S^1$. Let 
$\beta\in \tilde{K}(S^2)$ be the Bott element and 
$J : \tilde{K}^{-1}(S(G/S^+))\to \pi^0_S(S(G/S^+))$ 
be the stable complex $J$-homomorphism. Let 
$h : S^{k+1}(G/S^+)=S^k(S(G/S^+)) \to S^k$ denote the map representing 
$J(\beta\xi)\in \pi^0_S(S(G/S^+))$ where $\xi$ represents its own 
stable equivalence class. Then we know from ~\cite[Lemmas 2.2, 2.3]{O}
that the composite  
\begin{equation}
g=-h\circ Sf : S^{k+n}=S(S^{k+n-1})\to S(S^k(G/S^+))\approx 
S^k(S(G/S^+))\to S^k, 
\end{equation}
represents $[G, L]\in \pi^S_n$, i.e.      
$[G, L]=-\langle J(\beta\xi), [G/S]\rangle$. 

Under the identification of $G/S$ above, the total spaces of $\nu$ 
can be regarded as being written $\mathbb{R}^k\times ((E_+\times H/S)
\cup_{\bar{\tau}}(E_-\times H/S))\subset \mathbb{R}^{n-1+k}$ 
and so $g$ can be written in the form  
\begin{equation}
g : S^{k+n}\xrightarrow{-Sf} S^{k+1}((E_+\times H/S)^+
\cup_{\bar{\tau}^+}(E_-\times H/S)^+) \xrightarrow{h}S^k. 
\end{equation} 

For convenience we write $S(\mathbb{R}^k\times G/S)$ instead of 
$S^{k+n}$, the domain of the map $g$ above. Then   
$S(\mathbb{R}^k\times (E_\pm\times H/S))\subset S^{k+n}$ 
can be taken to be the hemispheres of $S^{k+n}$ by considering as
$\mathbb{R}^k\times (E_\pm\times H/S)\subset 
\mathbb{R}^k\times(\mathbb{R}_\pm^{m+1}\times\mathbb{R}^{n-m-2})
=\mathbb{R}_\pm^{n+k-1}$ where $\mathbb{R}_\pm^{m+1}$ are the half spaces of 
$\mathbb{R}^{m+1}$ consisting of $(t_0, \cdots, t_m)\in\mathbb{R}^{m+1}$ with 
$\pm t_m\ge 0$ (the sign $\pm$ applies in the same order as written). 
\begin{proof}[Proof of Theorem 1] 
Let us put $S_\pm^m=E_\pm/\partial E_\pm$.
Then due to the assumption that $\pi_{m-1}^S=0$ we find that 
the expression of $g$ in (2) can be rewritten as  
\[
g : S^{k+n}\xrightarrow{-(Sf)'}S^{k+1}((E_+\times H/S)^+
\cup_{\bar{\tau}'^+} (S_-^m\times H/S)^+) \xrightarrow{h'}S^k  
\]
where $-'$ denotes the map induced by $-$ and further in particular 
$\bar{\tau}'$ is the map of $\partial E_+\times H/S$ to 
$b_-\times H/S\subset\partial E_-\times H/S$ given by 
$\bar{\tau}'(x, yS)=(b_-, T(x)yS)$, $b_\pm$ denoting the collapsed 
$\partial E_\pm$. 
Here in above, the subgroup $K\subset H$ contains a subgroup isomorphic to $SU(2)$, so identifying it with $SU(2)$ we take $S$ to be $U(1)\subset SU(2)$ where $SU(2)/S\approx S^2$. Now in order to observe the behavior of 
$\bar{\tau}'$ we replace $H/S$ above by $SU(2)/S$ and regard $\bar{\tau}'$ 
as a map from $\partial E_+\times SU(2)/S$ to $b_-\times SU(2)/S$. This 
makes sense because of $\pi_2(H/S)\cong\mathbb{Z}$ which follows from the 
exact sequence of homotopy groups for the fibering $H\to H/S$. Then 
if we let $g_\pm$ denote the restrictions of $g$ to $S(\mathbb{R}^k\times  (E_\pm\times H/S))\subset S^{n+k}$, then we see that the value of $g_-$ can be represented as $d(q+1)$ times the value of $g_+$ for some $q$, $d$ given above. This multiple number can be interpreted as meaning that $H/S(=b_+\times H/S)$ overlaps on $H/S(=b_-\times H/S)$ 
$d(q+1)$ times under the deformation of $g$ above; in particular, $q+1$ expresses the degree of overlap itself and $d$ indicates how many times it occurs. 

Now the assumption on $\pi_{m-1}^S$ asserts slightly more strongly that $g$  satisfies
\begin{equation}
g_-\!\mid\!S(\mathbb{R}^k\times  (\partial E_-\times H/S)) \simeq c_\infty
\end{equation}
with the notation above where $c_\infty$ denotes the constant map at the base point. Applying this we see that the map g above can be further deformed into the composite
\begin{equation} 
g : S^{n+k}\xrightarrow{-(Sf)''} S^{k+1}(S_+^m\wedge H/S^+)
\vee S^{k+1}(S_-^m\wedge H/S^+)
\xrightarrow{h''_+\vee h''_-} S^k\vee S^k\xrightarrow{\mu} S^k
\end{equation}
where $-''$ also denotes the map induced by $-'$, $\mu$ the folding map; here from the relation bewteen $g_\pm$ observed above we have   
\begin{equation}
h''_-\simeq (-1)^{m-1}d(q+1)h''_+.
\end{equation}

Furthermore consider replacing these $h''_\pm$ by the maps
\[ \tilde{h}''_\pm : S^{k+1}(S_\pm^m\wedge H/S^+)\to S_\pm^m\wedge S^k, 
\quad (t, x, yS)\to x\wedge h''_\pm(t, x, yS) .\]
Then the composition (4) is transformed into the form 
\begin{equation} 
\begin{split}
g' : S^{n+k}&\xrightarrow{-(Sf)''} S^{k+1}(S_+^m\wedge H/S^+)
\vee S^{k+1}(S_-^m\wedge H/S^+)\\
&\approx S^{m+k+1}(H/S^+)\vee S^{m+k+1}(H/S^+)
\xrightarrow{\tilde{h}''_+\vee \tilde{h}''_-} S^{m+k}\vee S^{m+k}
\xrightarrow{S^m\mu} S^{m+k} 
\end{split} 
\end{equation}
through a canonical homeomorphism. From the construction we see that 
$\tilde{h}''_\pm$ become homotopic to the $m$-fold suspension of maps 
$h_\pm : S^{k+1}(H/S^+) \to S^k$ where each of them represents the map corresponding to $h$ in $(*)$ with $G$ replaced by $H$. From (5) it also follows 
that they must satisfy 
$h_-\simeq (-1)^{m-1}d(q+1)h_+$. By definition, applying the condition that 
$[G, L]=0$, i.e., $g\simeq c_\infty$ in the construction of $g'$, we find 
$g'\simeq c_\infty$. Thus we have  
 $(1+(-1)^{m-1}d(q+1))[H, L]=0\in \pi_{n-m}^S$. This proves 
the theorem.  
\end{proof}

\section{Proof of Theorem 2}
Note that we use the calculation results of the homotopy groups of spheres 
due to ~\cite{T} and ~\cite{MT} without reference. 

\begin{proof}
i) $G=SU(2)$. Since $G/S=S^2$, $\tilde{K}^{-1}(S(G/S))\cong \mathbb{Z}\beta(\xi-1)$; then by ~\cite{A} we have
\[\pi_S^0(S^3)=\mathbb{Z}_{24}J(\beta(\xi-1)).\] 
Consider the standard embedding of $G/S=S^2$ into $\mathbb{R}^3$. 
Then in a similar way as above, via the Pontrjagin-Thom construction, we 
have a stable map $f : S^{2+k}\to S^k(G/S^+)=S^k({S^2}^+)$ such that its homotopy class represents  $[G/S]$. If we take $f$ to be that in (1), then 
$\langle J(\beta(\xi-1), [G/S]\rangle=-[G, L]$. But in the present case,  
due to the construction of $[G/S]$,  we have that $\langle J(\beta(\xi-1), [G/S]\rangle$ must be identical to $J(\beta(\xi-1))$ and therefore 
\[J(\beta(\xi-1))=-[G, L].\] 
This together with the above equation tells us that 
$\pi^S_3=\mathbb{Z}_{24}[G, L]$.

ii) $G=SU(3)$.  Take $H=SU(2)$ and $K=S=U(1)$. Then $G/H=S^5$  
and $\pi_4^S=0$. From ~\cite[\S 24.3]{St} we know that 
$p\circ T : S^4\to H\to H/K=S^2$ is essential and so $d=1$ since 
$\pi_4(S^2)=\mathbb{Z}_2$. 
Now $2[G, L]=0\in \pi_8^S=\mathbb{Z}_2\oplus\mathbb{Z}_2$ and therefore 
applying Theorem 1 with $[G, L]$ replaced by twice itself we have 
$2(q+2)[H, L]=0$ where $1\le q\le 11$ since $\pi_6(S^2)=\mathbb{Z}_{12}$. 
But since ${\rm ord}([H, L])=24$ by i) above it follows that $q+2$ must be  
divisible by 12, so $q$ becomes equal to 10; hence it follows that
\[d=1, \ q=10.\] 

Suppose $[G, L]=0$. Then substituting these values into (1)  we have 
$12[H, L]=0$ which implies that the order of $[H, L]$ is reduced by at least half. This is clearly a contradiction. Hence we must have $[G, L]\ne 0$ which shows that $\pi_8^S=\mathbb{Z}_2
\oplus\mathbb{Z}_2[G, L]$.

iii) $G=Sp(2)$. Take $H=Sp(1)=SU(2)$ and put $K=S=U(1)$. 
Then $G/H=S^7$ and by ~\cite[\S\S24.3, 24.5]{St} we know that 
$p\circ T : S^6\to H \to H/K=S^2$ represents a nonzero element of 
$\pi_6(S^2)=\mathbb{Z}_{12}$, so we have $1\le d \le 11$. 
Let $T' : S^6\to SU(3)$ be 
the characteristic map of the bundle $SU(4)\to S^7$. Then by 
~\cite[\S24.5]{St} we see that $i\circ T \simeq T'$ where 
$i : H\hookrightarrow SU(3)$ is the inclusion and by ~\cite[\S25.2]{St}, 
using the results of ~\cite{Ker}, we also see that $T'$ represents 
a generator of $\pi_6(SU(3))=\mathbb{Z}_6$ because of $\pi_6(SU(4))=0$. 
From these two facts it follows that $T$ is twice a generator of 
$\pi_6(H)=\mathbb{Z}_{12}$, so that the value of $d$ can be reduced to 
\[1\le d \le 5.\] 

Now since $\pi_6^S=\mathbb{Z}_2$, $2\pi_6^S=0$. This equation permits us to
apply Theorem 1 with $[G, L]$ replaced by six times itself to 
$6[G, L]=0\in \pi_{10}^S=\mathbb{Z}_6$. Then by $(*)$ we have    
$6(1+d(q+1))[H, L]=0$ where $q=0$ or 1 since $\pi_8(S^2)=\mathbb{Z}_2$. 
Here by i) above ${\rm ord}([H, L])=24$ and so it is clear that 
$1+d(q+1)$ must be divisible by 4, which makes it possible to obtain 
\[d=3, \ q=0.\]

Similarly, if we suppose $2[G, L]=0$, then we have $2(1+d(q+1))[H, L]=0$ 
under the same condition as above, i.e. under the condition that $2\pi_6^S=0$, 
so by substituting the values of $d$, $q$ obtained above into $(*)$ in the case 
where $[G, L]$ replaced by twice itself we have $8[H, L]=0$. This contradicts the fact that ${\rm ord}[H, L]=24$; therefore we have $2[G, L]\ne 0$.  

Next consider the non-zeroness of $3[G, L]$. For this we observe 
$g : S^{10+k}\to S^k$ represented by (2). Let $S^5=S^6\cap S^6_\perp$ 
where $S^6_\perp\subset S^7$ denotes the equator defined by $t_0=0$, 
i.e. $S^5$ consists of the elements of $(0, t_1, \cdots, t_6, 0)\in S^7$. 
Then its restriction to 
$S(\mathbb{R}^k\times S^5 \times H/S)\subset S^{10+k}$ 
becomes null homotopic since $\pi_5^S=0$. 
Therefore we see that $g\!\mid\!S(\mathbb{R}^k\times S^6_\bot \times H/S)$ 
does homotopic to the sum of twice a map, due to the symmetry 
property of $T$ in the $t_0, t_1, \cdots, t_6$ coordinates. 
But since $\pi_9^S=\mathbb{Z}_2\oplus\mathbb{Z}_2\oplus\mathbb{Z}_2$ 
it also becomes null homotopic. By use of this, due to the symmetry 
property of $T$ again, $g$ can be written as $g\simeq 2a$ where 
$a : S^{10+k}\to S^k$. Hence $3g\simeq 6a\simeq c_\infty$ 
because of $\pi_{10}^S=\mathbb{Z}_6$. 
This shows that $3[G, L]=0$. Combining this with the result that $2[G, L]\ne 0$ just obtained above we can conclude that 
$\pi_{10}^S=\mathbb{Z}_2\oplus\mathbb{Z}_3 [G, L]$.

iv) $G=G_2$. There is a subgroup $H=SU(3)$ such that $G/H=S^6$. Here 
$\pi_5^S=0$ which shows that the given condition is satisfied. 
From ~\cite{Y} we recall that there is an inclusion homomorphism 
$G_2\to SO(7)$ such that $K$ and $H$, where $K=SU(2)$, are mapped into 
$SO(5)$ and $SO(6)$ as subgroups, respectively, keeping their inclusion relations 
$K\subset H\subset G_2$ and $SO(5)\subset SO(6)\subset SO(7)$. Then we also 
know that $T : S^5 \to H$ becomes homotopic in $SO(6)$ to the characteristic map $T' : S^5\to SO(6)$ of the bundle $SO(7)\to S^6$. By  ~\cite[\S23.4]{St}, 
$p'\circ T' : S^5\to S^5$ has degree 2 where $p': SO(6) \to S^5$ 
is the quotient map, so $p\circ T : S^5 \to S^5$ has also degree 2. 
Therefore $d$ must be a multiple of 2. 

Suppose now that $[G, L]=0$. Then since $d$ is even, substituting it into (1) 
we have $[H, L]=0$ since ${\rm ord}([H, L])=2$ by ii) above. This is 
a clear contradiction. So we must have $[G, L]\ne 0$ and therefore we 
can conclude that $\pi_{14}^S=\mathbb{Z}_2\oplus\mathbb{Z}_2 [G, L]$.

v) $G=SU(4)$. Take $H=SU(3)$. Then $G/H=S^7$. Since $\pi_5^S=0$, 
taking into account the symmmetry property in the coordinates 
$t_0, \cdots, t_6$, we see that the restriction of $g : S^{15+k}\to S^k$ to 
$S(\mathbb{R}^k\times \partial E_- \times H/S)\subset S^{15+k}$ 
becomes homotopic to the sum of twice a map. Moreover since 
$\pi_6^S=\mathbb{Z}_2$, this restriction map becomes null-homotopic. 
This shows that in the present case, the triviality of $\pi_6^S$, i.e. the 
given condition does not satisfied, but the equation (3) is satisfied. 
From this, in view of the proof of Theorem 1, we see that $(*)$ is 
applicable to $g$ above.

From ~\cite[\S24.3]{St} we know that $p\circ T : S^6\to H\to H/K=S^5$, 
where $K=SU(2)$, is inessential, so $d=0$. Hence assuming $[G, L]=0$ 
we have $[H, L]=0$ from $(*)$. This contradicts the fact that 
${\rm ord}([H, L])=2$ in ii) above, so it must be that $[G, L]\ne0$. 

From the observation in  the proof of the case iv) above we see that
$T\!\mid\!S^5$ can be taken to be the characteristic map $T'$ of the 
bundle $G_2\to S^6$ by looking at the matrix form of $T$ under 
$\pi_5(SO(6))=\mathbb{Z}$  ~\cite[p.131]{St}. So both of 
$p\circ T' : S^5\to H \to S^5$ and 
$(p\circ T)\!\mid\! S^5 : S^5\to H \to S^5$ have degree 2. Since 
$\pi_5(SU(3))\cong\pi_5(U)=\mathbb{Z}$ it follows that $T\!\mid\!S^5\simeq T'$. Taking this fact with the symmetry property of the restriction map of $g$ observed above we have $[G, L]=[S\times G_2, L]$ which shows that  
$[G, L]=\eta[G_2, L]$; therefore we can conclude that   
$\pi_{15}^S=\mathbb{Z}_{420}\oplus\mathbb{Z}_2[G, L]$.  

vi) $G=Sp(3)$. Take $H=Sp(2)$ and $K=Sp(1)$. Then $G/H=S^{11}$. In 
a similar way as in the case v) above we see that $3g$ satisfies the equation (3) above, using the equations $\pi_5^S=0$, $\pi_6^S=\mathbb{Z}_2$, 
$\pi_8^S=\mathbb{Z}_2\oplus \mathbb{Z}_2$ and $\pi_{10}^S=\mathbb{Z}_6$ 
and taking account into the symmetry property of $g_- : S^{21 +k}\to S^k$.  
Hence by (4) we obtain a decomposition such that $3g\simeq a_+ + a_-$ where 
$a_\pm : S^{21+k}\to S^k$ and $a_-\simeq d(q+1)a_+$. It therefore follows that 
\[g\simeq (1+d(q+1))a_+\] since $\pi_{21}^S=\mathbb{Z}_2\oplus\mathbb{Z}_2$.
From the same reasoning we find that the argument for (6) is applicable to $g$. 
Suppose $[G, L]=0$, i.e. $g\simeq c_\infty$. Then we have 
$(1+d(q+1))[H, L]=0$ by $(*)$. Since $\mathrm{ord}[H, L]=3$ by iii) above 
it follows that $1+d(q+1)$ is a multiple of 3. By 
~\cite[\S24.5]{St} and ~\cite[\S24.3]{St} we know that $d=1$ and by 
$\pi_{12}(S^2)=\mathbb{Z}_2\oplus\mathbb{Z}_2$ we also know that 
$0\le q\le 3$. Hence it must be that $d=q=1$. Substituting these values into 
the above equation we see that $g\simeq c_\infty$ is equivalent to 
$a_+\simeq c_\infty$. In the same way applying the argument for (6) to 
the latter equation there we have $[H, L]=0$. This is clearly a contradiction.
Hence we see that $g$ is not null homotopic and so 
$\pi_{21}^S=\mathbb{Z}_2\oplus\mathbb{Z}_2 [G, L]$.
\end{proof}

\end{document}